\documentclass[12pt]{article}
\usepackage{graphics, latexsym}
\usepackage{graphicx}
\usepackage{setspace}
\usepackage{amssymb}
\usepackage{amsmath}
\usepackage{epsfig}
\usepackage{CJK}

\begin{document}
\title{ Convergence Rate of  Sample  Mean for $\varphi$-Mixing Random Variables with Heavy-Tailed Distributions$^*$}
\author{TANG Fuquan $\,\,\,\,\,\,$  HAN Dong$^*$ \\
Department of Statistics, School of Mathematical Sciences, \\
Shanghai Jiao Tong University, Shanghai, 200240, China  }

\maketitle
\begin{abstract}
This article studies the convergence rate of the sample mean for $\varphi$-mixing dependent random variables with finite means and infinite variances. Dividing the sample mean into sum of the average of the main parts and the average of the tailed parts, we not only obtain the convergence rate of the sample mean but also prove that the convergence rate of the average of the main parts is faster than that of the average of the tailed parts.
\end{abstract}
\renewcommand{\thefootnote}{\fnsymbol{footnote}}
\footnotetext{
${}^*$Supported by National Natural Science Foundation of China (11531001)
\newline ${}^*$ Corresponding author, E-mail: donghan@sjtu.edu.cn }



\section{Introduction}
 A kind of phenomenon which often happens in  many random systems is that the mean is finite but the variance infinite. For example, consider the log return of Dow Jones Industrial  Average   $X_k=100\log (Y_k/Y_{k-1})$, where $Y_k$  denotes the closed price at $k$th day.  By using the data from 2019/04/26 to 2020/02/11, we can check that they are not independent and  the tailed index of the Hill estimator can be estimated as $2.1292$. This means that the mean is finite but the variance infinite.  One of questions which we are concerned here is that what is the convergence rate of the sample mean  $n^{-1}\sum_{k=1}^nX_k$ in probability when $X_k, 1\leq k\leq n,$ are mutually dependent and all the variances are infinite ?

Let us recall some known results. Let $X_n, n\geq 1,$  be i.i.d. sequence with the finite $1+\alpha$ moment $\textbf{E}|X_k|^{1+\alpha}=\mu_{\alpha}$ for some $\alpha \in (0, 1]$. Let $\delta \in (0, 1)$. Bubeck et al.${}^{[2]}$ proposed a truncated empirical mean $\hat{\mu}_T=n^{-1}\sum_{k=1}^nX_kI(|X_k|\leq b_k)$ to estimate the mean $\mu$ and proved that the convergence rate of $|\hat{\mu}_T-\mu|$ is $O(\log(1/\delta)/n)^{\alpha/(1+\alpha)}$ with the probability at least $1-\delta$, where $b_k=(bk/\log(1/\delta))^{1/(1+\alpha)}$, $b\geq \mu_{\alpha}$ and $I(.)$ is the indicator function. Lee et al.${}^{[6]}$ further presented a novel robust estimator with the same convergence rate, where information about $\mu_{\alpha}$ is not required in prior.  Avella-medina et al.${}^{[1]}$ considered Huber's M-estimator${}^{[4]}$ $\hat{\mu}_H$ of $\mu$ which is defined as the solution to $\sum_{k=1}^n\psi_H(X_k-\mu)=0$, where $\psi_H(x)=\min\{H, \max(-H, x)\}$ is the Huber function. They also obtained the same convergence rate of $|\hat{\mu}_H-\mu|$ when $H=(\mu_{\alpha}n/\log(2\delta^{-1}))^{1/(1+\alpha)}$.

We know that the tailed part of a random variable $X$ with the heavy-tailed distribution can best reflect the characteristics of probability. However, the above authors  did not consider the convergence rate of  the average of the tailed parts $n^{-1}\sum_{k=1}^n[X_kI(|X_k|>b_n)-\textbf{E}(X_kI(|X_k|>b_n))]$, where $b_n \to \infty$ as $n\to \infty$. This is a problem worth studying.  Another problem said by Bubeck et al${}^{[2]}$ is that  "It is unclear whether similar results could be obtained for heavy-tailed bandits with dependent reward processes", that is, the known results above heavily rely on the independence of the processes.

In this paper, we will investigate the convergence rate of the sample mean for $\varphi$-mixing dependent random variables with heavy-tailed distributions and try to solve the above two problems by a different method.

\section{Convergence rate of the sample mean }
Let  $\{X_n, n\geq 1\}$ be a random sequence satisfying  the following $\varphi-$mixing condition
\begin{eqnarray}
\varphi(m)&\leq & C (m+1)^{-r}  \text{\, for some constant $r>2$ and a positive constant $C$. } \\
\varphi(m)&=&\max_{k\geq 1, \,\, A \in \mathfrak{F}_1^{k},\,\,\textbf{P}(A)>0,\,\, B\in \mathfrak{F}_{k+m}^{\infty}}|\textbf{P}(B|A)-\textbf{P}(B)|\nonumber
\end{eqnarray}
for $m\geq 1$, where $\varphi(m)$ is decreasing on $m$, $\textbf{P}(\cdot)$ denotes the probability measure for the sequence $\{X_{k}, \, k\geq 0\}$ and $\mathfrak{F}_k^{m}=\sigma(X_{i}, \, k\leq i\leq m)$ for $1\leq k\leq m\leq +\infty$. From (1) it follows that
\begin{eqnarray*}
\varphi: =\sum_{m=1}^{\infty}\varphi^{1/2}(m)<\infty.
\end{eqnarray*}
Assume that there is a positive constant $ \alpha $ satisfying  $0<\alpha<1$ such that
\begin{eqnarray}
0<\nu_{\alpha}=\max_{n\geq 1}\textbf{E}(|X_n|^{1+\alpha})<\infty.
\end{eqnarray}

\textbf{Remark 1}.  If $\{X_n, n\geq 1\}$  has  a common marginal heavy-tailed distribution with $0<\textbf{E}(|X_n|^{1+\alpha})<\infty$, then the condition (2)  holds.

In order to get the convergence rate of the average of the tailed parts,  we first extend the result of the von Bahr-Esseen moment inequality for pairwise independent random variables in Theorem 2.1 in  Chen et al.${}^{[3]}$ to the case of $\varphi-$mixing random variables.

\textbf{Lemma 1.} \textit{  Let the conditions (1)-(2) hold. There is a positive constant $c_{\alpha}$ not depending  on $n$ such that
\begin{eqnarray}
\textbf{E}|\sum_{k=1}^nZ_k(b)|^{1+\alpha}\leq 2^{1+\alpha}c_{\alpha}\nu_{\alpha} n.
\end{eqnarray}}

\textbf{Proof.}  Let  $Z_k=Z_k(b)$,  $s=1+\alpha$ and  $M_n=\sum_{k=1}^n\textbf{E}(|Z_k|^{s})$. For any $\varepsilon >0$, by the same method as proving (2.1)-(2.3) in Theorem 2.1 in Chen et al.${}^{[3]}$, we can get that
\begin{eqnarray}
\textbf{E}|\sum_{k=1}^nZ_k|^{s}&\leq &(2+\varepsilon)M_n+J_n
\end{eqnarray}
and
\begin{eqnarray*}
\sup_{x\geq (1+\varepsilon)M_n}x^{-1/s}|\sum_{k=1}^n\textbf{E}[Z_kI(|Z_k|\leq x^{1/s})]|&=& \sup_{x\geq (1+\varepsilon)M_n} x^{-1/s}|\sum_{k=1}^n\textbf{E}[Z_kI(|Z_k|> x^{1/s})]|\\
&\leq &\sup_{x\geq (1+\varepsilon)M_n} x^{-1/s} x^{1/s-1}\sum_{k=1}^n\textbf{E}(|Z_k|^{s}I(|Z_k|> x^{1/s}))\\
&\leq & (1+\varepsilon)^{-1},
\end{eqnarray*}
where
\begin{eqnarray*}
J_n=\int_{(1+\varepsilon)M_n}^{\infty}\textbf{P}\Big( |\sum_{k=1}^n[Z_kI(|Z_k|\leq x^{1/s})]|>x^{1/s}\Big)dx.
\end{eqnarray*}
Let $\tilde{\varepsilon}=[1-(1+\varepsilon)^{-1}]^{-2}$. From Markov's inequality, (1.11) in Wang et al.${}^{[7]}$ and the inequality  of $\varphi$-mixing (Ibragimov${}^{[5]}$), $|\textbf{Cov}(\xi, \eta)|\leq 2[\varphi(m)Var(\xi)Var(\eta)]^{1/2}$, it follows that
\begin{eqnarray}
J_n &\leq &\int_{(1+\varepsilon)M_n}^{\infty}\textbf{P}\Big( |\sum_{k=1}^n[Z_kI(|Z_k|\leq x^{1/s})-\textbf{E}(Z_kI(|Z_k|\leq x^{1/s}))]|>[1-(1+\varepsilon)^{-1}]x^{1/s}\Big)dx\nonumber\\
&\leq & \tilde{\varepsilon}\sum_{k=1}^n\int_{(1+\varepsilon)M_n}^{\infty}x^{-2/s}\textbf{E}\Big(Z^2_kI(|Z_k|\leq x^{1/s})\Big)dx\nonumber\\
&& +  2\tilde{\varepsilon}\sum_{k=1}^n\int_{(1+\varepsilon)M_n}^{\infty}x^{-2/s}\Big(\sum_{j=k+1}^n\textbf{Cov}\Big(Z_kI(|Z_k|\leq x^{1/s}), Z_jI(|Z_j|\leq x^{1/s})\Big)\Big)dx\nonumber\\
&&\leq \tilde{\varepsilon}(1+4\varphi )\sum_{k=1}^n\int_{(1+\varepsilon)M_n}^{\infty}x^{-2/s}\textbf{E}\Big(Z^2_kI(|Z_k|\leq x^{1/s})\Big)dx.
\end{eqnarray}
Furthermore, like proving (2.5)-(2.7) in Theorem 2.1 in Chen et al.${}^{[3]}$, we have
\begin{eqnarray}
\int_{(1+\varepsilon)M_n}^{\infty}x^{-2/s}\textbf{E}\Big(Z^2_kI(|Z_k|\leq x^{1/s})\Big)dx\leq \Big(\frac{2s}{(2-s)^2}+\frac{2}{2-s}\Big)\textbf{E}(|Z_k|^{s}).
\end{eqnarray}
Hence, from (4) (5) and (6) it follows that
\begin{eqnarray}
\textbf{E}|\sum_{k=1}^nZ_k|^{1+\alpha}&\leq &\Big(2+\varepsilon+\frac{1+4\varphi }{[1-(1+\varepsilon)^{-1}]^{2}}\times \frac{4}{(1-\alpha)^2}\Big)M_n \nonumber\\
&=& f(\varepsilon)M_n.
\end{eqnarray}
Note that the function $f(\varepsilon)$ in (7) is positive and continuous,  and $f(\varepsilon)$ goes  to $\infty$ as $\varepsilon \searrow 0$ and $\varepsilon\nearrow \infty$, respectively. This means that there is a positive constant $c_{\alpha}$ such that $c_{\alpha}=\min_{\varepsilon}\{f(\varepsilon)\}$. Thus, the inequality (3) follows from (7) since $\textbf{E}(|Z_k(b)|^{1+\alpha})\leq 2^{1+\alpha}\nu_{\alpha}$ for all $k\geq 1$. It completes the proof.

\textbf{Remark 2}. It is clear that when  $X_n, n\geq 1,$  are mutually independent and $\textbf{E}(X_n^{2})<\infty$, the inequality (3) holds for $\alpha =1$, where $c_{\alpha}=1/4$.

Let $Y_k(b)=X_kI(|X_k|\leq b)-\textbf{E}(X_kI(|X_k|\leq b))$ for $b>0$ and $k\geq 1$. Here $n^{-1}\sum_{k=1}^nY_k(b)$ can be considered as the average of the main parts. Let
\begin{eqnarray*}
\beta=\frac{\alpha r +1}{(1+\alpha)(1+r)},\,\,\,\,\,\,\,\, \gamma= \frac{\alpha r +\alpha'}{(1+\alpha)(1+r)}, \,\,\,\,\,\,\,\,\, A=e^{1/4+\sqrt{e}C4^r},
\end{eqnarray*}
where $0\leq \alpha'\leq \alpha$.

The following theorem shows the convergence rate of the sample mean in probability.

\textbf{Theorem 1.} \textit{  Assume that the conditions (1)-(2) hold. Let $c$ be a positive constant,  $\delta \in (0, 1)$, $p=2\log (1/\delta)$,  $a_1=\sqrt{3c\nu_{\alpha}(1+4\varphi )}$ and $a_2=2(c_{\alpha}\nu_{\alpha})^{1/(1+\alpha)}$. Then, there is a positive number $n_0$ depending on $\delta$ such that
\begin{eqnarray}
\textbf{P}\Big(|\frac{1}{n}\sum_{k=1}^n[X_k -\textbf{E}(X_k)]| \geq a_1(\frac{p}{n})^{\beta} +a_2(\frac{p}{n})^{\gamma}\Big)\leq 2A\delta^c+p^{-\alpha}\Big(\frac{p}{n}\Big)^{\frac{\alpha-\alpha'}{1+r}}
\end{eqnarray}
for $n\geq n_0$.}

\textbf{Proof.}  Let $b_n=(n/p)^{\lambda}$, $Y_k=Y_k(b_n)=X_kI(|X_k|\leq b_n)-\textbf{E}(X_kI(|X_k|\leq b_n))$ and $Z_k=Z_k(b_n)=X_kI(|X_k|>b_n)-\textbf{E}(X_kI(|X_k|>b_n))$, where $\lambda=(r-1)[(1+r)(1+\alpha)]^{-1}$. Dividing $X_k -\textbf{E}(X_k)$ into two parts: the main part $Y_k$ and the tailed part $Z_k$, we have
\begin{eqnarray}
&&\textbf{P}\Big(|\frac{1}{n}\sum_{k=1}^n[X_k -\textbf{E}(X_k)]| \geq a_1(\frac{p}{n})^{\beta} +a_2(\frac{p}{n})^{\gamma}\Big)\nonumber\\
&&\leq \textbf{P}\Big(|\frac{1}{n}\sum_{k=1}^nY_k| \geq a_1(\frac{p}{n})^{\beta}\Big)+ \textbf{P}\Big(|\frac{1}{n}\sum_{k=1}^nZ_k|\geq a_2(\frac{p}{n})^{\gamma}\Big)
\end{eqnarray}
By  Chebyshev's inequality and Lemma 1 we have
\begin{eqnarray}
\textbf{P}\Big(|\frac{1}{n}\sum_{k=1}^nZ_k|\geq a_2(\frac{p}{n})^{\gamma}\Big)\leq \frac{\textbf{E}|\sum_{k=1}^nZ_k|^{1+\alpha}}{(a_2n)^{1+\alpha}(p/n)^{(1+\alpha)\gamma}}\leq p^{-\alpha}\Big(\frac{p}{n}\Big)^{\frac{\alpha-\alpha'}{1+r}},
\end{eqnarray}
where $O((p/n)^{\gamma})$ is the convergence rate of the average of the tailed part in probability.

Let $S_n=\sum_{k=1}^nY_k$. As in Ibragimov ${}^{[5]}$ and Yang ${}^{[8]}$, the sum $S_n$ can be written as two parts
\begin{eqnarray}
S_n =\sum_{i=1}^{l+1}\xi_i+\sum_{i=1}^{l}\eta_i,\,\,\,\,\,\,\,\,\,\,\xi_i=\sum_{j=1}^mY_{2(i-1)m+j},\,\,\,\,\,\,\,\,\eta_i= \sum_{j=1}^mY_{(2i-1)m+j}
\end{eqnarray}
for $1\leq i\leq l$, $\xi_{l+1}=0$ for $2lm\geq n$ and $\xi_{l+1}=Y_{2lm+1}+...+Y_{n}$ for $2lm < n$, where
\begin{eqnarray}
m=\lceil\frac{1}{M^2}(\frac{n}{p})^{1/(r+1)}\rceil,\,\,\,\,\,\,\,\, l=\lceil M^2\log (1/\delta)(\frac{n}{p})^{r/(r+1)}\rceil,
\end{eqnarray}
where  $\lceil x \rceil$ denotes the maximum integer of $x$ and $M$ is a large positive number satisfying that for given a small $0<\epsilon \leq 1/2$, $1/M \leq \epsilon$. It is clear that $2lm=2\lceil n/2\rceil\approx n$ for large $n$.

Let $\delta_n=a_1(p/n)^{\beta}$. From Chebyshev's inequality and H\"{o}lder inequality it follows that
\begin{eqnarray}
&&\textbf{P}\Big(\sum_{k=1}^nY_k \geq n\delta_n\Big) \leq  \exp\{-\theta n\delta_n \}\textbf{E}\exp\{\theta S_n\}\nonumber \\
&&\leq  \exp\{-\theta n\delta_n \}\Big(\textbf{E}\exp\{2\theta \sum_{i=1}^{l+1}\xi_i\}\Big)^{1/2}\Big(\textbf{E}\exp\{2\theta \sum_{i=1}^{l}\eta_i\}\Big)^{1/2}
\end{eqnarray}
for $\theta >0$.  Furthermore, by the inequality  of $\varphi$-mixing and (1.11) in Wang et al.${}^{[7]}$,  we have
\begin{eqnarray*}
\textbf{E}(\xi^2_i)&=& \sum_{j=1}^m\textbf{E}(Y^2_{2(i-1)m+j})+2\sum_{j=1}^m\sum_{j'=j+1}^m\textbf{E}(Y_{2(i-1)m+j}Y_{2(i-1)m+j'})\\
&\leq & (1+4\varphi )\sum_{j=1}^m\textbf{E}(X^2_{2(i-1)m+j}I(|X_{2(i-1)m+j}|\leq b_n))\\
&=&  (1+4\varphi )\sum_{j=1}^m\textbf{E}(|X_{2(i-1)m+j}|^{1+\alpha}|X_{2(i-1)m+j}|^{1-\alpha}I(|X_{2(i-1)m+j}|\leq b_n))\\
&\leq & m\nu_{\alpha}(1+4\varphi )(b_n)^{1-\alpha}.
\end{eqnarray*}
Note that $2\theta |\xi_i|\leq 2\theta m b_n\leq 1/M \leq \epsilon\leq 1/2 $ for $1\leq i\leq l+1$ when $\theta \leq M(p/n)^\kappa$, where $\kappa=(\alpha +r)[(1+r)(1+\alpha)]^{-1}$. Let $\theta \leq M(p/n)^\kappa$. Hence,
\begin{eqnarray*}
\log \textbf{E}(\exp\{2\theta \xi_i\})=\log \Big(1+ 2\theta^2\textbf{E}(\xi^2_i)(1+\epsilon)\Big)\leq 2\theta^2\textbf{E}(\xi^2_i)(1+\epsilon)
\end{eqnarray*}
and therefore,
\begin{eqnarray}
\textbf{E}(\exp\{2\theta \xi_i\})&\leq &\exp\{2\theta^2\textbf{E}(\xi^2_i)(1+\epsilon)\} \nonumber\\
&\leq & \exp\{2\theta^2m\nu_{\alpha}(1+4\varphi )(b_n)^{1-\alpha}(1+\epsilon)\}.
\end{eqnarray}
By using (14), the definition of $\{\xi_i\}$ in (11), the $\varphi-$mixing property and the similar method used by Yang ${}^{[8]}$,  we can get that
\begin{eqnarray}\label{eq31}
&&\textbf{E}(\exp\{2\theta \sum_{i=1}^l \xi_i\})= \textbf{E}\Big(\exp\{2\theta \sum_{i=1}^{l-1}\xi_i\}\textbf{E}( \exp\{2\theta \xi_l\}|\mathfrak{F}_0^{l-1})\Big)\nonumber\\
&\leq & (\textbf{E}( \exp\{2\theta \xi_l\})+e^{\epsilon}\varphi(m))\textbf{E}(\exp\{2\theta \sum_{i=1}^{l-1}\xi_i\})\nonumber\\
&\leq & \exp\{e^{\epsilon}\varphi(m)l\}\exp\{2\theta^2ml\nu_{\alpha}(1+4\varphi C_2)(b_n)^{1-\alpha}(1+\epsilon)\}.
\end{eqnarray}
Similarly, we can obtain that
\begin{eqnarray}
\textbf{E}(\exp\{2\theta \sum_{i=1}^l \eta_i\})\leq \exp\{e^{\epsilon}\varphi(m)l\}\exp\{2\theta^2ml\nu_{\alpha}(1+4\varphi )(b_n)^{1-\alpha}(1+\epsilon)\}.
\end{eqnarray}
Note that $2\theta |\xi_{l+1}|\leq 1/M \leq \epsilon$ for $\theta \leq M(p/n)^\kappa$, $2lm\approx n$, $\varphi(m)l\leq C  m^{-r}l\approx CM^{2r}$ and
\begin{eqnarray*}
 e^{\epsilon/2}\exp\{e^{\epsilon}\varphi(m)l\}\leq A
\end{eqnarray*}
for $\epsilon=1/2$ and $M=2$. It follows from (13), (15)-(16) that
\begin{eqnarray*}
&&\textbf{P}\Big(\sum_{k=1}^nY_k \geq n\delta_n\Big) \leq  e^{\epsilon/2}\exp\{e^{\epsilon}\varphi(m)l\}\exp\{-\theta n\delta_n  \}\exp\{2\theta^2ml\nu_{\alpha}(1+4\varphi )(b_n)^{1-\alpha}(1+\epsilon)\}\nonumber\\
&\approx & e^{\epsilon/2}\exp\{e^{\epsilon}\varphi(m)l\}\exp\{-n[\delta_n\theta -(1+\epsilon)\nu_{\alpha}(1+4\varphi )(b_n)^{1-\alpha}\theta^2]\}\nonumber\\
&\leq &A \exp\{-n[\delta_n\theta -(1+1/2)\nu_{\alpha}(1+4\varphi)(b_n)^{1-\alpha}\theta^2])\}.
\end{eqnarray*}
Note that the function $g(\theta)=\delta_n\theta -(1+1/2)\nu_{\alpha}(1+4\varphi )(b_n)^{1-\alpha}\theta^2$  arrives its maximum value at $\theta^*_n=\delta_n[3\nu_{\alpha}(1+4\varphi )(b_n)^{1-\alpha}]^{-1}$,
$\delta^2_n(b_n)^{-(1-\alpha)}=6c\nu_{\alpha}(1+4\varphi )\log (1/\delta)/n$ and $\kappa<1$. Take the minimum number $n_0$ such that $\theta^*_{n_0}\leq 2(p/n_{0})^\kappa$, we can get that
\begin{eqnarray}
\textbf{P}\Big(\sum_{k=1}^nY_k \geq n\delta_n\Big) &\leq &  A\exp\{-n(\delta^2_n[6\nu_{\alpha}(1+4\varphi )(b_n)^{1-\alpha}]^{-1}\}\nonumber\\
     &\leq &  A\exp\{-c\log(1/\delta)=A\delta^c
\end{eqnarray}
for $n\geq n_0$. Similarly, we have
\begin{eqnarray}
\textbf{P}\Big(\sum_{k=1}^nY_k <- n\delta_n\Big)\leq A\delta^c.
\end{eqnarray}
This means that the convergence rate of the average of the main parts in probability is $O((p/n)^{\beta})$.  By (9), (10), (17) and (18), we see that  the inequality (8) holds. It completes the proof of Theorem 1.

\textbf{Remark 3}. It can be seen from (8), (9) and (10) that the convergence rate of the average of the main parts  is faster than that of the average of the tailed parts since $\beta >\alpha/(1+\alpha)\geq \gamma$ for $r<\infty$. Furthermore, taking $\alpha'=\alpha$, we see that the  convergence rate of the average of the tailed parts can arrive at $O((p/n)^{\alpha/(1+\alpha)})$ with the probability at least $1-p^{-\alpha}$.

In the following three corollaries we assume that  $\mu=\textbf{E}(X_k)$ for all $k\geq 1$.
Let $\alpha'<\alpha$. Take $c_0$ and $n_1$ such that $2A\delta^{c_0}=\delta/2$ and $p^{-\alpha}(p/n_1)^{(\alpha-\alpha')/(1+r)}=\delta/2$. We have
the following corollary.

\textbf{Corollary 1.} \textit{ Let  $\alpha'<\alpha$ and $c\geq c_0$. Then, with the probability at least $1-\delta$, one has
\begin{eqnarray*}
|\frac{1}{n}\sum_{k=1}^nX_k -\mu| < a_1[\frac{2\log (1/\delta)}{n}]^{\beta} +a_2[\frac{2\log (1/\delta)}{n}]^{\gamma}
\end{eqnarray*}
for $n\geq \max\{n_0,\, n_1\}$.}

Note that $\gamma=\alpha/(1+\alpha)$ when $\alpha'=\alpha$. Hence we have the following lemma.

\textbf{Corollary 2.} \textit{ Let   $\alpha'=\alpha$ and $c\geq c_0$. Then, with the probability at least $1-\delta/2-(2\log (1/\delta))^{-\alpha}$, one has
\begin{eqnarray*}
|\frac{1}{n}\sum_{k=1}^nX_k -\mu| < a_1[\frac{2\log (1/\delta)}{n}]^{\beta} +a_2[\frac{2\log (1/\delta)}{n}]^{\alpha/(1+\alpha)}
\end{eqnarray*}
for $n\geq n_0$.}

It is known that  $X_n, n\geq 1,$ will become mutually independent when $r=\infty$ in (1).

\textbf{Corollary 3.} \textit{ If $X_n, n\geq 1,$ are mutually independent, that is, $r=\infty$ in (1) and therefore, $a_1=\sqrt{3c\nu_{\alpha}}$ and $\beta=\gamma=\alpha/(1+\alpha)$, then, with the probability at least $1-\delta/2-(2\log (1/\delta))^{-\alpha},$ one has
\begin{eqnarray*}
|\frac{1}{n}\sum_{k=1}^nX_k -\mu| < (a_1+a_2)\Big(\frac{2\log (1/\delta)}{n}\Big)^{\alpha/(1+\alpha)}
\end{eqnarray*}
for $n\geq n_0$ and $c\geq c_0$.}

\textbf{Remark 4.} If $X_n, n\geq 1,$ are i.i.d. sequence with $\textbf{E}(|X_n|^{1+\alpha})<\infty$ for $0\leq \alpha \leq 1$, then the conditions (1)-(2) hold.

\section{ Conclusion}

By the different method we obtain  the convergence rate of the sample mean for $\varphi$-mixing dependent random variables with finite means and infinite variances. Furthermore, we prove that the convergence rate $O((2\log (\delta^{-1})/n)^{\beta})$ of the average of the main parts is faster than the convergence rate  $O((2\log (\delta^{-1})/n)^{\alpha/(1+\alpha)})$ at which  the average of the tailed parts can arrive.

If $X_n, n\geq 1,$ are i.i.d. sequence with $\textbf{E}(|X_n|^{1+\alpha})<\infty$ for $0< \alpha \leq 1$, it follows from the corollary 3 and Remark 3 that we not only get the same results obtained by Bubeck et al.${}^{[2]}$ and Avella-medina et al.${}^{[1]}$, but also show that the convergence rate of the average of the tailed parts is $O((2\log (\delta^{-1})/n)^{\alpha/(1+\alpha)})$.

\par
\par

\textbf{Acknowledgements}   We sincerely thank two reviewers for their precious comments on the manuscript.

$\,\,\,\,\,\,\,\,\,$\\
$\,\,\,\,\,\,\,\,\,$\\
$\,\,\,\,\,\,\,\,\,$\\

%
%
%
%
%
%
%
%
%

\end{document}